\documentclass[12pt]{amsart}
\usepackage{epsfig}
\usepackage{amsmath}
\usepackage{amsfonts}
\usepackage{amssymb}
\usepackage{epigraph, fancyhdr}

\def\F{\mathcal F}
\def\K{\mathcal K}
\def\J{\mathcal J}
\def\L{\mathcal L} 

\def\O{\mathcal O}

\def\1{\mathbf 1}
\def\M{{\overline{\mathcal M}}}

\def\QQ{\mathbb Q}
\def\ZZ{\mathbb Z}
\def\CC{\mathbb C}

\def\Res{\operatorname{Res}}

\def\tilde{\widetilde}

\def\a{\alpha}

\def\t{{\mathbf t}}

\def\gl{\nu}
\def\gL{\Lambda}

\def\lan{\langle}
\def\ran{\rangle}

\def\ev{\operatorname{ev}}
\def\ft{\operatorname{ft}}

\def\td{\operatorname{td}}
\def\ch{\operatorname{ch}}

\def\tr{\operatorname{tr}}

\renewcommand{\Delta}{\triangle}

\title[Elementary K-theory of $\M_{0,n}/S_n$]
      {Permutation-equivariant \\ quantum K-theory I. \\
      Definitions \\ Elementary K-theory of $\M_{0,n}/S_n$}

\author[A. Givental]{Alexander GIVENTAL}
\thanks{This material is based upon work supported by the National 
Science Foundation under Grant DMS-1007164, and by the IBS Center for Geometry 
and Physics, POSTECH, Korea.} 

\date{June 23, 2015}

\begin{document}

\begin{abstract} 

K-theoretic Gromov-Witten (GW) invariants of a compact K\"ahler manifold $X$ are defined as super-dimensions of sheaf cohomology of interesting bundles over moduli spaces of $n$-pointed holomorphic curves in $X$. In this paper, we introduce K-theoretic GW-invariants cognizant of the $S_n$-module structure on the sheaf cohomology, induced by renumbering of the marked points, and compute some of these invariants for $X=pt$.   
  
\end{abstract}

\maketitle

\section*{Preface}

In Fall 2014, I gave a talk on the subject of permutation-equivariant quantum K-theory and its relations to mirror symmetry at {\em The Legacy of Vladimir Arnold} conference in Toronto. After the talk, explaining to Anatoly Vershik why the work was not published yet, I got a piece of good advice from him: he suggested that one should publish not a whole theory, but little pieces of it.

The present paper begins a series of such pieces. Each one is supposed to have its own punch-line, and be reasonably self-contained, or at least readable separately from the others. Yet, they are chapters of the same story, follow a single plan, and are meant {\em to be continued}.
One of our intentions is to identify the right place for toric $q$-hypergeometric functions among genus-0 K-theoretic Gromov--Witten invariants.  Another one is to elucidate the role of finite-difference operators. In particular, we will see that the $q$-exponential function is even more prominent in quantum K-theory then the ordinary exponential function is in quantum cohomology. As a remote goal, we would like the $q$-analogues of the Witten--Kontsevich tau-function to arise from K-theory of the Deligne--Mumford quotients $\M_{g,n}/S_n$. 
 
This and several forthcoming chapters are based on the lectures I gave in June-July 2015 at the IBS Center for Geometry and Physics at POSTECH, Korea. I'd like to thank the Center's director Yong-Geun Oh and his staff for their hospitality and for creating an ideal working environment.   

\section*{$S_n$-equivariant correlators}

Let $X$ be a compact K\"ahler manifold, a {\em target space} of GW-theory,
$X_{g,n,d}$ denote the moduli space of degree-$d$ stable maps to $X$ of nodal
compact connected $n$-pointed curves of arithmetical genus $g$, 
$\ev_i:X_{g,n,d} \to X$ the evaluation map at the $i$th marked point, $L_i$
the line bundle over $X_{g,n,d}$ formed by the cotangent lines to the curves at
the $i$th marked point. Given elements $\phi_i \in K^0(X)$ and integers 
$k_i\in\ZZ$, $i=1,\dots, n$, one defines a K-theoretic GW-invariant of $X$
as the holomorphic Euler characteristic
\[  \lan \phi_1 L^{k_1},\dots, \phi_n L^{k_n}\ran_{g,n,d} := 
\chi \textstyle \left (X_{g,n,d} ; \O_{g,n,d}^{virt}\otimes \prod_{i=1}^n 
L_i^{k_i} \ev_i^*(\phi_i)\right).\]
Here $\O_{g,n,d}^{virt}$ is the {\em virtual structure sheaf} 
introduced by Y.-P. Lee 
\cite{YPLee} as the K-theoretic counterpart of virtual fundamental cycles in the
cohomological theory of GW-invariants. The above ``correlators'' can be extended
poly-linearly to the space of Laurent polynomials 
\[ t(q)=\sum_{m\in \ZZ} t_m q^m, \ \ t_m=\sum_{\a} t_{m,\a}\phi_{\a}\]
(here $\{ \phi_{\a}\}$ is a basis in $K^0(X)\otimes \QQ$, and $t_{k,\a}$ are
formal 
variables),  and thereby encode the values of all individual correlators
by the totally symmetric degree-$n$ polynomial
$\lan t(L), \dots, t(L)\ran_{g,n,d}$.          

Our aim is to enrich this information using the action of $S_n$ by 
permutations of the marked points. Namely, since the marked points are numbered,
their renumbering on a given stable map produces a new stable map, and hence
this operation induces an automorphism of the moduli space:
$X_{g,n,d}\to X_{g,n,d}$. 
In fact the automorphism is relative over $X_{g,0,d}$ (here we have in mind the
map
$\ft: X_{g,n,d}\to X_{g,0,d}$ defined by  forgetting the marked point). The map 
$\ft$ respects the construction \cite{YPLee} of virtual structure sheaves: 
\[ \O_{g,n,d}^{virt}=\ft^*\O^{virt}_{g,0,d}.\]
Therefore, as long as the {\em inputs} $t(q)$ in all seats of the correlator
are the same,
the corresponding sheaf cohomology, and hence their alternated 
sum, carries a well-defined structure of a virtual $S_n$-module.
Let us introduce for this $S_n$-module the notation 
\[ [ t(L), \dots, t(L) ]_{g,n,d} := \sum (-1)^m H^m(X_{g,n,d}; \O_{g,n,d}^{virt}
\otimes \prod_{i=1}^n (\sum_{k\in \ZZ}\ev_i^*(t_k)L_i^k)).\]

Thus defined GW-invariants with values in the representation ring of $S_n$ lack 
two features required by the standard combinatorial framework of GW-theory: they are not poly-linear, and they take incomparable values for different values of $n$. We handle both difficulties by employing Schur--Weyl reciprocity.

Let $\gL$ be a $\lambda$-algebra, by which we will understand an algebra over $\QQ$ equipped with abstract {\em Adams operations} $\Psi^m$, $m=1,2,\dots $, i.e. ring homoorphisms $\gL\to \gL$ satisfying $\Psi^r \Psi^s = \Psi^{rs}$ and
$\Psi^1=\operatorname{id}$. \footnote{One usually defines
  $\lambda$-algebras in terms of axiomatic exterior power operation. For us the Adams operations will be more important. The difference, disappears over $\QQ$. The reason is that Newton polynomials are expressed as polynomials with integer coefficients in terms of elementary symmetric functions, but the inverse formulas involve fractions.} The following construction of correlators has direct topological meaning when $\gL = K^0(Y)\otimes \QQ$, the K-ring of some space $Y$ equipped with the natural Adams operations, but it can be extended
to arbitrary $\lambda$-algebras.

On the role of {\em inputs} we take Laurent polynomials
$\t =\sum_{m\in \ZZ} \t_m q^m$ with vector coefficients $\t_k \in K^0(X)\otimes \gL$. Given several such inputs $\t^{(1)},\dots \t^{(s)}$, we define correlators of {\em permutation-equivariant} quantum K-theory with several groups of sizes $k_1+\dots+k_s=n$ of identical inputs (and hence symmetric with respect to
the subgroup $H=S_{k_1}\times \cdots \times S_{k_s}$ of $S_n$), and taking values
in $\gL$:
\begin{align*}
  \lan &\t^{(1)},\dots, \t^{(1)}; \dots ; \t^{(s)},\dots,\t^{(s)}\ran_{g,n,d}^H := \\ &\left(\pi: (X_{g,n,d}\times Y) / H
  \to Y\right)_*\left(\O_{g,n,d}^{virt} \otimes \prod_{a=1}^s \prod_{i=1}^{k_a}
  \left( \sum_{m\in \ZZ}\ev_i^*(\t^{(a)}_m)L_i^m \right) \right)
  \end{align*}
where $\pi_*$ is the K-theoretic push-forward along the indicated projection map $\pi$. Note that the sheaf on the right lives naturally on $X_{g,n,d}\times Y$ and is $H$-invariant. Taking the quotient, by definition, extracts $H$-invariants from the K-theoretic push-forward to $Y$.

\medskip

{\tt Example 1.} {\em $GL_N$-equivariant K-theory.} Take $\gL$ to be the algebra of symmetric functions in $N$ variables $x_1,\dots, x_N$ with the Adams operations $\Psi^r(x_i)=x_i^r$. 
It can be viewed as (a subring in) $Repr \,GL_N = K^0(BGL_N(\CC))$, the representation ring of $GL_N(\CC)$, by considering $x_i$ to be the eigenvalues of diagonal
matrices in the vector representation $\CC^N$. Respectively, $K^0(X)\otimes \gL$ can be interpreted as $GL_N$-equivariant K-ring of $X$ equipped with the trivial $GL_N$-action. Let $t$ be a legitimate input of the ordinary quantum K-theory, i.e. Laurent polynomial $L$ with coefficients
from $K^0(X)$, and $\gl\in \gL$. Then 
\[ \lan \gl t,\dots, \gl t\ran_{g,n,d}^{S_N} =
\frac{1}{n!} \sum_{h\in S_n} \tr_h [t,\dots,t]_{g,n,d}\prod_{r=1}^{\infty}
  \Psi^r(\gl)^{l_r(h)},\]
where $l_r(h)$ denotes the number of cycles of length $r$ in the permutation $h$. Indeed, if $\gl$ in the correlator stands for a $GL_N$-module attached at each marked point, then $[ \gl t,\dots, \gl t]_{g,n,d}=[t, \dots t]_{g,n,d}\otimes \gl^{\otimes n}$. The second factor here is a $GL_N(\CC) \times S_n$-module. For a diagonal matrix $x$ and a permutation $h$, we have
\[ \tr_{(x,h)} \gl^{\otimes n} =\tr_x \prod_{r=1}^{\infty}\Psi^r(\gl)^{l_r(h)}.\]
Indeed, due to the universality of Adams operations, it suffices to check this for $\gl=\CC^N$, the vector representation, which is straightforward:
\[ \tr_{(x,h)} (\CC^N)^{\otimes n} = \prod_{r=1}^{\infty}N_r^{l_r(h)}(x),\]
where $N_r(x) =x_1^r+\cdots+x_N^r=\Psi^r(N_1)$ is the $r$th Newton polynomial. 

\medskip

{\tt Example 2:} {\em Schur--Weyl's reciprocity.}
According to Schur--Weyl's reciprocity, the $GL_N\times S_n$-character of 
$(\CC^N)^{\otimes n}$ has the form:
\[ \prod_{r=1}^{\infty} N_r^{l_r(h)}(x)
= \sum_{\Delta} w_{\Delta}(h) s_{\Delta}(x),\]
where $s_{\Delta}$ is the Schur polynomial, the character of the irreducible
$GL_N$-module with the highest weight determined by the partition 
(or the Young diagram) $\Delta$, and $w_{\Delta}$ is the character of the irreducible $S_n$-module corresponding to the same Young diagram. The diagrams here consist of $n$ cells and have no more than $N$ rows. The Schur polynomials $s_{\Delta}$
form a real orthonormal basis in the space of all symmetric polynomials 
of degree $n$ (in $N$ variables). Therefore, using the notation 
$(\cdot, \cdot)$ for pairing of representations (or characters), we have:  
\[ ( \lan t(L),\dots, t(L)\ran^{S_n}_{g,n,d}, s_{\nabla} ) = 
\frac{1}{n!}\sum_{h\in S_n}
\tr_h [t(L),\dots,t(L)]_{g,n,d} w_{\nabla}(h),\] 
that is, equal to the multiplicity of the irreducible $S_n$-module $\nabla$ in
the $S_n$-module of our interest. 

\medskip

{\tt Example 3:} {\em $N\to \infty$.} In this limit, $\gL$ becomes the abstract algebra of symmetric functions $\QQ[[ N_1, N_2,\dots ]]$ with
the Adams operations $\Psi^r(N_m)=N_{rm}$. This example captures the entire
information about $[t,\dots,t]_{g,n,d}$ as $S_n$-modules for all $n$ simultaneously.

\medskip

{\tt Example 4:} {\em Symmetrized quantum K-theory.} Taking in Example 1 $N=1$,
we obtain $\gL=\QQ [x]$ with $\Psi^r(x)=x^r$. This choice corresponds extracting $S_n$ invariants from sheaf cohomology:
\[ [xt,\dots, xt]_{g,n,d} = [t,\dots, t]_{g,n,d}^{S_n} x^n.\]
Indeed, the action of $S_n$ on $GL_1$-module $(\CC^1)^{\otimes n}$ is trivial.
We will refer to this important special case of permutation-equivariant
quantum K-theory as {\em permutation-invariant} or {\em symmetrized}.

\medskip

{\tt Example 5:} {\em The permutation-equivariant binomial formula.} Returning to the definition of permutation-equivariant correlators, we can see that they 
possess permutation-equivariant version of poly-additivity. For instance,
\[ \lan \t'+\t'',\dots, \t'+\t''\ran^{S_n}_{g,n,d} = \sum_{k+l=n}
\lan \t',\dots \t', \t'',\dots \t''\ran^{S_k\times S_l}_{g,n,d}.\]
Using the bracket notation $\lan \dots \ran$ for the sheaf cohomology on $X_{g,n,d}\times Y$ (i.e. before taking $S_n$-invariants), we have
the following equality of $S_n$-modules:
\[ \lan \t'+\t'',\dots, \t'+\t''\ran_{g,n,d} = \sum_{k+l=n} \operatorname{Ind}_{S_k\times S_l}^{S_n} \lan\t',\dots,\t',\t'',\dots, \t''\ran_{g,n,d},\]
where $\operatorname{Ind}_H^G$ denotes the operation of inducing a $G$-module
from an $H$-module. Extracting $S_n$-invariants on both sides proves the claim.
Indeed, due to the reciprocity between inducing and restricting, for any $H$-module $V$, we have $(\operatorname{Ind}_H^GV)^G = V^H$, since restricting the trivial $G$-module to $H$ yields the trivial $H$-module. 

\medskip

Finally introduce the genus-$g$ {\em descendent potentials} of
permutation-equivariant quantum K-theory:
\[ \F_g = \sum_{d,n} Q^d \lan \t(L),\dots,\t(L)\ran_{g,n,d}^{S_n}.\] 
Here $Q^d$ is, as usual, the monomial representing the degree $d \in H_2(X)$ in
the Novikov ring. Note that the customary in Taylor's formulas division by $n!$
is replaced by extracting $S_n$-invariants. The potential is a formal function
on the space of Laurent polynomials in $q$ with coefficients in $K^0(X)\otimes \gL$. We assume that $\lambda$-algebra $\gL$ is extended to power series in Novikov's variables (e.g. one could take $\gL=\QQ [[N_1,N_2,\dots ]] \, [[Q]]$) and the Adams operations are extended by $\Psi^r(Q^d)=Q^{rd}$. We will refer to $\gL$ as {\em Newton-Novikov's ring}.  

{\tt Remark.} I am thankful to A. Polishchuk, who pointed out to me that in a related context of {\em modular operads}, an equivalent formalism of encoding permutation-equivarint information using the algebra of symmetric functions \cite{M} was used by E. Getzler and M. Kapranov \cite{GeK}.    

\section*{The small J-function of the point}

In this section, we use an explicit description of Deligne--Mumford spaces 
$\M_{0,n}$ in terms of Veronese curves to compute the ``small''  J-function 
in the permutation-equivariant quantum K-theory of $X=pt$.

\medskip

{\tt Theorem.} {\em For  $\gl \in \gL$, put 
\[ J_{pt}(\gl):=1-q+\gl+\sum_{n\geq 2} 
\lan  \gl,\dots,\gl, \frac{1}{1-qL}\ran^{S_n}_{0,n+1}.\]
Then
\[ J_{pt}=(1-q)e^{\sum_{k>0}\Psi^k(\gl)/k(1-q^k)}.\]}

\medskip

{\tt Proof.} We refer to the paper \cite{Kap} by M. Kapranov for details of
the description of $\M_{0,n+1}$ in terms of Veronese curves in $\CC P^{n-2}$,
i.e. generic rational curves of degree equal to the dimension of the
ambient projective space. They are all isomorphic to the model Veronese curve 
$(u:v)\mapsto (u^{n-2}:u^{n-3}v:\dots : uv^{n-3}:v^{n-2})$ under the action of 
$PGL_2(\CC)\times 
PGL_{n-1}(\CC)$ by reparameterizations and projective automorphisms, and form 
a family of dimension $(n+1)(n-3)$. The moduli space $\M_{0,n+1}$ is identified 
with a suitable closure of the space of Veronese curves passing through a fixed
generic configuration of $n$ points $p_1,\dots, p_n\in \CC P^{n-2}$. According
to \cite{Kap}, the closure can be taken in the Chow scheme of algebraic cycles
(or in the suitable Hilbert scheme). Moreover, $\M_{0,n+1}$ is obtained
explicitly by a certain succession of blow-ups of $\CC P^{n-2}$ centered at all 
subspaces passing through the $n$ points. The rational map, inverse to the 
projection $\pi: \M_{0,n+1}\to \CC P^{n-2}$, can be
described this way: for a generic $p\in \CC P^{n-2}$, there is a unique Veronese
curve passing through $(p_1,\dots,p_n,p)$. (Example: a unique conic through $5$ generic points on the plane.)

The forgetting map $\ft_{n+2}: \M_{0,n+2}\to \M_{0,n+1}$ can be described as follows
(see Figure 1).
Veronese curves of degree $n-1$ in $\CC P^{n-1}$ passing through fixed generic points $p_1,\dots, p_n, p_{n+1}$ can be projected from $p_{n+1}$ to $\CC P^{n-1}$, to become Veronese curves of degree $n-2$ passing through the projections
$\tilde{p}_1,\dots, \tilde{p}_n$ of $p_1,\dots, p_n$. 
According to \cite{Kap}, this projection survives the passage to the Chow closure. 

\begin{figure}[htb]
\begin{center}
\epsfig{file=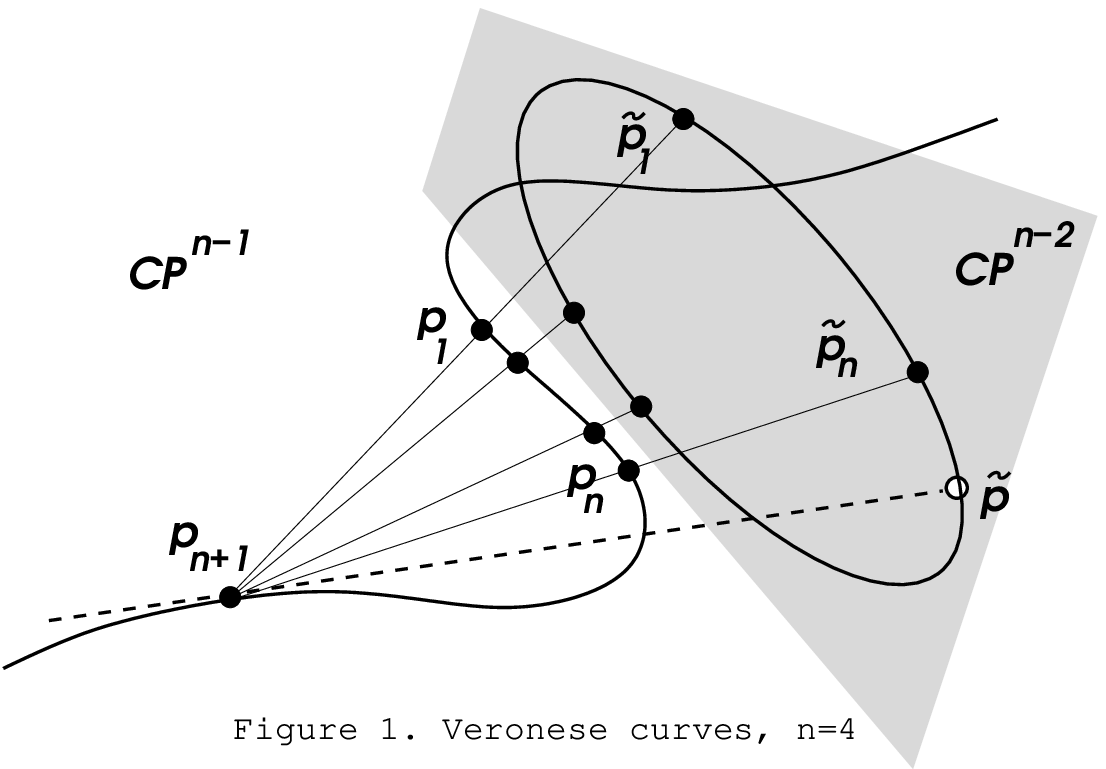}
\end{center}
\end{figure}

Moreover, as it follows from the exact description of the succession of the blow-ups (see \cite{Kap}, Theorem 4.3.3), the section
$\M_{0,n+1}\subset \M_{0,n+2}$ of the forgetful map $\ft_{n+2}: \M_{0,N+2}\to \M_{0,n+1}$ defined by the $n+1$-st marked point is obtained by blowing up $\CC P^{n-1}$ at
$p_{n+1}$, and then taking the proper transform of the exceptional divisor $\CC P^{n-2}$ under all further blow-ups. Their centers come from higher-dimensional subspaces passing through $p_{n+1}$, and are transverse to the the divisor. This means that conormal bundle to the section (which is the official definition of $L_{n+1}$ over $\M_{0,n+1}$) coincides with the pull-back of the conormal bundle to the exceptional $\CC P^{n-2}$ (which is $\O(1)$) by the blow-down map $\pi: \M_{0,n+1}\to \CC P^{n-2}$. Thus, $L_{n+1} = \pi^*\O(1)$.

Note that since $L_{n+1} = \pi^* \O(1)$, then $\pi_* L_{n+1}^m = \O(m)$, because 
the K-theoretic push-forward of the structure sheaf along a blow-down map has trivial higher direct images. Thus the problem of computing $J_{pt}$ receives the following elementary 
interpretation. Let $S_n$ act on $\CC P^{n-2}=\operatorname{proj}(\CC ^{n-1})$ by 
permutations of the vertices $p_1,\dots, p_n$ of the standard simplex. Then 
the $S_n$-module denoted in the previous section $[1,\dots, 1, L^m]_{0,n+1}$ 
is the space of degree $m$ polynomials in $\CC^{n-1}$. Respectively,
\[ \lan \gl,\dots \gl,\frac{1}{1-qL} \ran_{0,n+1}^{S_n} = 
\frac{1}{n!}\sum_{h\in S_n} \tr_h S^*_q(\CC ^{n-1}) \prod_{r>0}\Psi^r(\gl)^{l_k(h)},\]
where $S^*_q(\CC^{n-1})=\oplus_{m\geq 0}q^mS^m(\CC^{n-1})^*$ is the graded (and 
weighted by powers of $q$) algebra of polynomial functions on $\CC^{n-1}$.
       
The series $J_{pt}$, the total sum of the correlators over all $n$, can be computed by Lefschetz fixed point formula. In fact summation over all symmetric groups 
can be rewritten in terms of conjugacy classes. The action of $h\in S_n$ 
on $\CC^n$ (rather than $\CC^{n-1}$ decomposes into the direct product of 
elementary $k$-cycles $c_k$ acting on $\CC^k$ by the cyclic permutation of 
the coordinates. The trace $\tr_{c_k} S^*_q(\CC^k)$ can be computed as 
$\prod_{s=1}^k(1-e^{2\pi i s/k}q)^{-1}=(1-q^k)^{-1}$, since $e^{2\pi i s}$ are simple 
eigenvalues of $c_k$ on $\CC^k$.
Taking in account the size $n!/\prod_k l_k! k^{l_k}$ of the conjugacy class with $l_k$ cycles of length $k$, we conclude that 
\[ \sum_{n\geq 0}\frac{1}{n!}\sum_{h\in S_n} 
\tr_h S^*_q(\CC ^n) \prod_{k>0}\Psi^k(\gl)^{l_k(h)} = 
\sum_{l_1,L_2,\dots >0} \prod_{k>0}\frac{1}{l_k!}\left(\frac{\Psi^k(\gl)}
{k(1-q^k)}\right)^{l_k}.\]
The latter sum coincides with $e^{\textstyle\sum_{k>0} \Psi^k(\nu)/k(1-q^k)}$.
The extra factor $(1-q)$ in the theorem takes care of the excess (comparing to 
$\CC^{n-1}$) $1$-dimensional subspace in $\CC^n$ with the trivial action of $S^n$, because the Poincar{\'e} polynomial $\tr_{id} S^*_q(\CC) = 1/(1-q)$. $\square$
    
\medskip

{\tt Corollary 1.} {\em In the symmetrized theory, the value of the J-function
\[ J_{pt}^{sym}:=1-q+x+\sum_{n\geq 2} x^n 
\dim [\frac{1}{1-qL}, 1,\dots, 1]_{0,n+1}^{S_n} \]
is expressed in terms of the $q$-exponential function 
$e_q(y):=\sum_{n\geq 0} \frac{y^n}{[n]_q!}$:
\[ J^{sym}_{pt}=(1-q) e_q\left(\frac{x}{1-q}\right)=
\sum_{n\geq 0}\frac{x^n}{(1-q^2)\dots (1-q^n)}.\]}

{\tt Proof.} Taking in the theorem $\gL=\QQ [[x]]$ (i.e. choosing $GL_N$ to be $GL_1$), and sertting $\gl=x$,  we find 
\[ f(x):=(1-q)^{-1}J^{sym}_{pt} = e^{\sum_{k>0}x^k/k(1-q^k)}.\]
Note that $f$ satisfies the following finite-difference equation:
\[ f(x)-f(qx) = f(x) \left(1-e^{-\sum_{k>0} x^k/k}\right) = f(x) (1-(1-x))
=xf(x).\]
For $e_q$, we also have:
\[ e_q(\frac{x}{1-q})-e_q(\frac{qx}{1-q}) = \sum_{n\geq 0} \frac{x^n(1-q^n)}
{(1-q)(1-q^2)\dots (1-q^n)} = x e_q(\frac{x}{1-q}).\]
Since both are power series in $x$ with the free term $1$, they coincide. 
$\square$

\medskip

{\tt Corollary 2.} {\em When $\gL$ is the algebra of symmetric functions in
  $x_1,\dots, x_N$, and $\gl=t N_1$, where $t$ is a scalar, we have
\[ J_{pt}(tN_1) = (1-q) \prod_i e_q\left(\frac{x_i}{1-q}\right)^t.\]}

{\tt Proof.} Write $\Psi^k(tN_1)=
t(x_1^k+\cdots+x_N^k)$ for each $k$. $\square$

\enddocument